\documentclass[11pt]{article}
\usepackage[margin=1in]{geometry} 
\usepackage{amssymb,amsfonts,amsmath,bbm,mathrsfs,stmaryrd}
\usepackage{xcolor}
\usepackage{url}

\usepackage{enumerate}

\usepackage{hyperref}
\hypersetup{colorlinks,
             linkcolor=black!75!red,
             citecolor=blue,
             pdftitle={},
             pdfproducer={pdfLaTeX},
             pdfpagemode=None,
             bookmarksopen=true
             bookmarksnumbered=true}
           
\usepackage{tikz}
\usetikzlibrary{arrows,calc,decorations.pathreplacing,decorations.markings,intersections,shapes.geometric,through,fit,shapes.symbols,positioning,decorations.pathmorphing}

\usepackage{braket}

\usepackage[amsmath,thmmarks,hyperref]{ntheorem}
\usepackage{cleveref}

\creflabelformat{enumi}{#2(#1)#3}

\crefname{section}{Section}{Sections}
\crefformat{section}{#2Section~#1#3} 
\Crefformat{section}{#2Section~#1#3} 

\crefname{subsection}{\S}{\S\S}
\AtBeginDocument{%
  \crefformat{subsection}{#2\S#1#3}%
  \Crefformat{subsection}{#2\S#1#3}%
}

\theoremstyle{plain}

\newtheorem{lemma}{Lemma}[section]
\newtheorem{proposition}[lemma]{Proposition}
\newtheorem{corollary}[lemma]{Corollary}
\newtheorem{theorem}[lemma]{Theorem}

\theoremstyle{nonumberplain}

\theoremstyle{plain}
\theorembodyfont{\upshape}
\theoremsymbol{\ensuremath{\blacklozenge}}

\newtheorem{definition}[lemma]{Definition}
\newtheorem{example}[lemma]{Example}
\newtheorem{remark}[lemma]{Remark}

\crefname{definition}{definition}{definitions}
\crefformat{definition}{#2definition~#1#3} 
\Crefformat{definition}{#2Definition~#1#3} 

\crefname{ex}{example}{examples}
\crefformat{example}{#2example~#1#3} 
\Crefformat{example}{#2Example~#1#3} 

\crefname{remark}{remark}{remarks}
\crefformat{remark}{#2remark~#1#3} 
\Crefformat{remark}{#2Remark~#1#3} 

\crefname{convention}{convention}{conventions}
\crefformat{convention}{#2convention~#1#3} 
\Crefformat{convention}{#2Convention~#1#3}

\crefname{lemma}{lemma}{lemmas}
\crefformat{lemma}{#2lemma~#1#3} 
\Crefformat{lemma}{#2Lemma~#1#3} 

\crefname{proposition}{proposition}{propositions}
\crefformat{proposition}{#2proposition~#1#3} 
\Crefformat{proposition}{#2Proposition~#1#3} 

\crefname{corollary}{corollary}{corollaries}
\crefformat{corollary}{#2corollary~#1#3} 
\Crefformat{corollary}{#2Corollary~#1#3} 

\crefname{theorem}{theorem}{theorems}
\crefformat{theorem}{#2theorem~#1#3} 
\Crefformat{theorem}{#2Theorem~#1#3} 

\crefname{enumi}{}{}
\crefformat{enumi}{(#2#1#3)}
\Crefformat{enumi}{(#2#1#3)}

\crefname{assumption}{assumption}{Assumptions}
\crefformat{assumption}{#2assumption~#1#3} 
\Crefformat{assumption}{#2Assumption~#1#3} 

\crefname{question}{question}{Questions}
\crefformat{question}{#2question~#1#3} 
\Crefformat{question}{#2Question~#1#3} 

\crefname{problem}{problem}{Problems}
\crefformat{problem}{#2problem~#1#3} 
\Crefformat{problem}{#2Problem~#1#3}

\crefname{equation}{}{}
\crefformat{equation}{(#2#1#3)} 
\Crefformat{equation}{(#2#1#3)}

\numberwithin{equation}{section}

\theoremstyle{nonumberplain}
\theoremsymbol{\ensuremath{\blacksquare}}

\newtheorem{proof}{Proof}

\newcommand\cC{{\mathcal C}}
\newcommand\cD{{\mathcal D}}

\newcommand\cF{{\mathcal F}}

\newcommand\cL{{\mathcal L}}
\newcommand\cM{{\mathcal M}}

\DeclareMathOperator{\id}{id}

\newcommand{\qedhere}{\mbox{}\hfill\ensuremath{\blacksquare}}

\title{Higher dualizability and singly-generated Grothendieck categories}
\author{Alexandru Chirvasitu}

\begin{document}

\date{}

\newcommand{\Addresses}{{
  \bigskip
  \footnotesize

  \textsc{Department of Mathematics, University at Buffalo, Buffalo,
    NY 14260-2900, USA}\par\nopagebreak \textit{E-mail address}:
  \texttt{achirvas@buffalo.edu}

}}

\maketitle

\begin{abstract}
  Let $k$ be a field. We show that locally presentable, $k$-linear categories $\cC$ dualizable in the sense that the identity functor can be recovered as $\coprod_i x_i\otimes f_i$ for objects $x_i\in \cC$ and left adjoints $f_i$ from $\cC$ to $\mathrm{Vect}_k$  are products of copies of $\mathrm{Vect}_k$. This partially confirms a conjecture by Brandenburg, the author and T. Johnson-Freyd.

  Motivated by this, we also characterize the Grothendieck categories containing an object $x$ with the property that every object is a copower of $x$: they are precisely the categories of non-singular injective right modules over simple, regular, right self-injective rings of type I or III.
\end{abstract}

\noindent {\em Key words: locally presentable category; dualizable; abelian category; Grothendieck category; regular ring; self-injective ring; non-singular module; type}

\vspace{.5cm}

\noindent{MSC 2020: 18C35; 16D50; 18A30; 18A35}



\section*{Introduction}

The initial motivation for the present note stems from \cite{bcjf}, which studies the notion of {\it dualizability} for particularly well-behaved categories: a property capturing a higher-categorical version of the ``rigidity'' specific to finite-dimensional vector spaces or, more generally, projective finitely-generated modules. The concept is of interest at least in part due to the pervasiveness of dualities (of objects, morphisms, functors, etc.) in topological quantum field theory: see e.g. \cite[\S VII]{bd} or \cite[\S 2.3]{lur}, where {\it dualizable} and {\it fully dualizable} objects in a higher category are introduced.

To recall the main characters, briefly (with more background and references in the main text), fixing a field $k$ throughout:

\begin{definition}\label{def:init}
  A {\it 2-$k$-module} is a $k$-linear locally presentable category in the sense of \cite[Definition 1.17]{ar}.

  2-$k$-modules form a complete and cocomplete self-enriched 2-category ${}_{\mathrm{Vect}}\cM$ with $k$-linear left adjoints as 1-morphisms and $k$-linear natural transformations as 2-morphisms; the notation is meant to suggest the higher-categorical setting: 2-$k$-modules are modules over $\mathrm{Vect}$.

  The {\it dual} $\cC^*$ is the 2-$k$-module $\mathrm{Hom}(\cC,\mathrm{Vect})$.
  
  A 2-$k$-module $\cC$ is {\it 1-dualizable} (or just plain {\it dualizable}, for brevity) if the canonical functor
  \begin{equation}\label{eq:can}
    \cC\boxtimes \cC^*\to \mathrm{End}(\cC)
  \end{equation}
  is an equivalence.  
\end{definition}
(see \Cref{se.prel} for more on the categorical tensor product `$\boxtimes$').

There is a bit of a dearth of dualizable 2-$k$-modules {\it not} of the form
\begin{equation*}
  \text{linear functors }\Gamma\to \mathrm{Vect}
\end{equation*}
for small $k$-linear categories $\Gamma$, leading us to conjecture in passing in \cite[Remark 3.6]{bcjf} that they are all of this form. This paper was initially an attempt to confirm this in particularly simple cases, i.e. the simplest one would try (see \Cref{def:dlz,th:naive} below): 

\begin{theorem}\label{th:naive0}
  Suppose the 2-$k$-module $\cC$ admits sets
  \begin{itemize}
  \item $x_i\in \cC$ of objects and
  \item $f_i\in \cC^*$ of functors
  \end{itemize}
  such that the identity functor decomposes as $\coprod_i x_i\otimes f_i$. Then, $\cC$ is a coproduct of copies of the category $\mathrm{Vect}$.
\end{theorem}

\Cref{def:dlz} terms the categories in \Cref{th:naive0} {\it naively dualizable}, because they achieve dualizability in the most straightforward way imaginable. In particular, one might consider the case where the set $\{i\}$ of indices in \Cref{th:naive0} is a singleton: the identity functor is of the form $x\otimes f$ for an object $x\in \cC$ and a left adjoint functor $\cC\to \mathrm{Vect}$. In that case, the result specializes further (\Cref{cor:simp}):

\begin{corollary}\label{cor:simp0}
  If for the non-trivial 2-$k$-ring we have $\id\cong x\otimes f$ for an object $x$ and a functor $f\in \cC^*$, then $\cC\simeq \mathrm{Vect}$.
\end{corollary}

An earlier version of the paper dealt with \Cref{cor:simp0} directly (rather than as a consequence of \Cref{th:naive0}), by first noting that in every category $\cC$ as in \Cref{cor:simp0} every object is a coproduct (direct sum) of copies of $x$. This seems to be a peculiar property of some interest of its own, which motivates \Cref{se:sing}, which classifies {\it Grothendieck} categories with this property (\Cref{th:complete}): 

\begin{theorem}\label{th:complete0}
  For a Grothendieck category $\cC$ the following conditions are equivalent:
  \begin{enumerate}[(a)]
  \item There is an object $x\in \cC$ such that all objects of $\cC$ are copowers (i.e. coproducts of copies) of $x$. 
  \item $\cC$ is equivalent to the category $\cM_A^r$ of non-singular right injective modules over a simple, regular, right self-injective ring that is either
    \begin{itemize}
    \item a division ring, or
    \item of type III in the sense of \cite[Definition following Theorem 10.10]{good}.
    \end{itemize}
  \end{enumerate}
\end{theorem}

Recall (e.g. \cite[\S 2]{ah}) that {\it non-singular} modules are those containing no elements annihilated by essential ideals. Categories of non-singular injective modules feature prominently in the study of {\it regular} rings (\Cref{subse:reg}): the main structure theorem (\cite[Chapter XII, Theorem 1.3]{stns}, say) characterizing Grothendieck categories all of whose objects are projective identifies these with $\cM_A^r$ for regular self-injective $A$.

$\cM_A^r$ can be characterized in a number of other ways: it is also, for instance, the category of injectives over $A$ that are summands of direct products $A^I$ for various index sets $I$ (\cite[Chapter XII, Theorem 1.3]{stns}).

Note the contrast between \Cref{cor:simp0} and \Cref{th:complete0}: the dualizability requirement of the former appears to further ``rigidify'' the category, imposing the condition that objects be copowers of $x$ functorially. On the other hand, the second, type-III branch of \Cref{th:complete0} can be regarded as a loosening of that constraint, allowing for the less intuitive behavior characteristic of type-III rings: objects $x$ isomorphic to their own double copowers $x\oplus x$ \cite[Corollary 10.17]{good}, etc.

\subsection*{Acknowledgements}

This work was partially supported through NSF grant DMS-2001128.

I am grateful for enlightening comments from Martin Brandenburg (who suggested the original problem) and Theo Johnson-Freyd on the contents of \Cref{se:dlz} and Ken Goodearl on those of \Cref{se:sing}.

\section{Preliminaries}\label{se.prel}

All categories and functors are assumed $k$-linear for a fixed field $k$. We write
\begin{equation*}
  \mathrm{Vect} = {}_k\mathrm{Vect}
\end{equation*}
for the category of $k$-vector spaces.

We denote categories of modules by $\cM$, decorated with the base ring on the left or right, depending on whether they are left (or respectively right) modules. $\cM_A$, for instance, is the category of right $A$-modules over a ring $A$.

Inequality symbols between objects typically denote the relation of being a subobject, and $\le_{\oplus}$ means `is a direct summand of'.

\subsection{2-modules}\label{subse:2rng}

For much of the material below we refer to the more extensive \cite[\S 2]{bcjf} as well as \cite[Chapter 1]{ar}, which covers valuable background on locally presentable categories. Taking that for granted, we pause here only to recall that locally presentable categories are a particularly pleasant class of category to work with; a sample:

\begin{itemize}
\item (One version of) Freyd's adjoint functor theorem (\cite[\S V.8, Theorem 2]{macl}) holds for any such category $\cC$, in the sense that a cocontinuous functor with domain $\cC$ is automatically left adjoint (as follows for instance from the cited theorem and \cite[Remark 1.56]{ar});
\item They form a complete and cocomplete 2-category, with left adjoints as 1-morphisms and natural transformations as 2-morphisms (\cite[Proposition 2.1.11]{cjf}).
\end{itemize}

All of the above goes through in the $k$-linear setting, as recalled in \Cref{def:init}. We note (e.g. \cite[Proposition 2.1.11]{cjf} ) that the (co)product of a family of 2-$k$-modules $\cC_i$ is simply their product as categories.  

\begin{remark}
  Writing $\mathrm{hom}(\cC,\cD)$ for the category of 1-morphisms between 2-$k$-modules, ${}_{\mathrm{Vect}}\cM$ is symmetric monoidal with respect to a 2-categorical tensor product `$\boxtimes$' with the property that $\cC\boxtimes \cD$ receives a universal bifunctor $\cC\times \cD\to \cC\boxtimes \cD$ that is separately linear and cocontinuous.

  See for instance \cite[Lemma 2.7]{bcjf}, which in turn cites \cite[Exercise 1.1]{ar} and \cite[\S 6.5]{kel}.
\end{remark}

\begin{definition}\label{def:dlz}
  $\cC$ is {\it simply (1-)dualizable} if the identity functor $\id\in \mathrm{End}(\cC)$ is the image 
  \begin{equation}\label{eq:dec}
    x\boxtimes f\mapsto \id
  \end{equation}
  through \Cref{eq:can} of a simple tensor for some $x\in \cC$ and $f\in \cC^*$.

  $\cC$ is {\it naively (1-)dualizable} if the identity functor $\id\in \mathrm{End}(\cC)$ is the image
  \begin{equation}\label{eq:decmult}
    \coprod_{i\in I} x_i\boxtimes f_i\mapsto \id
  \end{equation}
  through \Cref{eq:can} of a coproduct of simple tensors for some family of $x_i\in \cC$ and $f_i\in \cC^*$.
\end{definition}

\subsection{Regular self-injective rings}\label{subse:reg}

We will have to dabble in the general theory of {\it von Neumann regular} rings, for which \cite{good} will serve as our main source. Recall (e.g. \cite[p.1, Definition]{good}) that these are the (always unital) rings $A$ with the property that for every $x\in A$ there is some $y\in A$ with $x = xyx$. Many equivalent characterizations exist:

\begin{itemize}
\item every one-sided (left or right) principal ideal is generated by an idempotent;
\item every one-sided (left or right) finitely generated ideal is generated by an idempotent;
\item every (left or right) module is flat (for which reason the rings are also termed {\it absolutely flat});
\item finitely generated submodules of projective (left or right) $A$-modules are direct summands. 
\end{itemize}

We compress the phrase `von Neumann regular' to just plain `regular' for brevity. Of special interest will be {\it right self-injective} regular rings, i.e. those regular $A$ which are injective as right modules over themselves. Unless specified otherwise, `self-injective' means `right self-injective'.

Regular right self-injective rings have received much attention in the literature; for our present purposes \cite[Chapters 9-12]{good} provide ample background, with \cite[Chapter 10]{good} of particular interest.

\section{Dualizability}\label{se:dlz}

\begin{theorem}\label{th:naive}
  A naively 1-dualizable 2-$k$-module $\cC$ is equivalent as a 2-$k$-module to a (co)product $\mathrm{Vect}^{\oplus S}$ of copies of $\mathrm{Vect}$.
\end{theorem}
\begin{proof}
  Suppose we have a decomposition \Cref{eq:decmult} of the identity functor on $\cC$. This then provides us with
  \begin{itemize}
  \item a functor $\pi:\mathrm{Vect}^{\oplus I}\to \cC$ defined by
    \begin{equation*}
      \mathrm{Vect}^{\oplus I}\ni (V_i)_{i\in I}\mapsto \bigoplus_{i\in I} x_i\otimes V_i \in \cC, and
    \end{equation*}
  \item a functor $\iota:\cC\to \mathrm{Vect}^{\oplus I}$ going in the opposite direction, defined as the product
    \begin{equation*}
      \iota:=\prod_i f_i:\cC\to \prod_{i\in I}\mathrm{Vect}\simeq \mathrm{Vect}^{\oplus I}.
    \end{equation*}
  \end{itemize}
  The dualizability assumption means precisely that $\iota$ splits $\pi$ in the sense that $\pi\circ\iota\cong \id$. 
  
  Since furthermore $\pi$ is a left adjoint (as are all morphisms of 2-$k$-modules), it must be right exact. On the other hand its domain $\mathrm{Vect}^{\oplus I}$ is a {\it spectral} Grothendieck category \cite[Chapter 4, Proposition 2.3]{pop}: all of its objects are projective (or equivalently, injective). This means that all short exact sequences split and hence $\pi$ is in fact exact (i.e. preserves both kernels and cokernels).

  Now denote by $\ker \pi$ the full subcategory of $\mathrm{Vect}^{\oplus I}$ consisting of objects annihilated by $\pi$ (cf. \cite[Chapitre III, Proposition 5]{gab} or \cite[\S 4.3, Exercise 5]{pop}). It is a {\it thick} or {\it dense} subcategory of $\mathrm{Vect}^{\oplus I}$ in the sense of \cite[\S III.1]{gab} (there called `\'epaisse') or \cite[\S 4.3]{pop}: given a short exact sequence
  \begin{equation*}
    0\to X'\to X\to X''\to 0,
  \end{equation*}
  $X$ belongs to $\ker\pi$ if and only if both $X'$ and $X''$ do.

  We can then form the {\it quotient} $\mathrm{Vect}^{\oplus I}/\ker\pi$ of $\mathrm{Vect}^{\oplus I}$ by $\ker\pi$ in the sense of \cite[discussion preceding Lemma 3.4]{pop}. Moreover, the canonical functor
  \begin{equation*}
    T:\mathrm{Vect}^{\oplus I} \to \mathrm{Vect}^{\oplus I}/\ker\pi
  \end{equation*}
  (of \cite[\S 4.3, Theorem 3.8]{pop}) is left adjoint, i.e. $\ker \pi$ is a {\it localizing} subcategory of $\mathrm{Vect}^{\oplus I}$ (\cite[p.372]{gab} or \cite[\S 4.4]{pop}); this follows, for instance, from \cite[\S III.4, Proposition 8]{gab}: $\ker\pi$ is clearly closed under taking colimits.
  
  According to \cite[\S 4.3, Exercise 5]{pop} we have a factorization
  \begin{equation}\label{eq:pipi}
    \begin{tikzpicture}[auto,baseline=(current  bounding  box.center)]
      \path[anchor=base] 
      (0,0) node (l) {$\mathrm{Vect}^{\oplus I}$}
      +(3,-.5) node (d) {$\mathrm{Vect}^{\oplus I}/\ker\pi$}
      +(6,0) node (r) {$\cC$}
      ;
      \draw[->] (l) to[bend left=6] node[pos=.5,auto] {$\scriptstyle\pi$} (r);
      \draw[->] (l) to[bend right=6] node[pos=.5,auto,swap] {$\scriptstyle T$} (d);
      \draw[->] (d) to[bend right=6] node[pos=.5,auto,swap] {$\scriptstyle \overline{\pi}$} (r);
    \end{tikzpicture}
  \end{equation}
  for some (unique, up to natural isomorphism) exact faithful functor $\overline{\pi}$.

  We do not quite know that $\overline{\pi}$ is an equivalence yet, because we do not, at the moment, have fullness; in fact, $\overline{\pi}$ will {\it not}, in general, be an equivalence: see \Cref{ex:notiso} below. It is, however,
  \begin{itemize}
  \item essentially surjective on objects, because it is right-invertible: $\overline{\pi}\circ T\circ \iota\cong \id$;
  \item faithful by construction.
  \end{itemize}
  According to \cite[\S III.4, Proposition 8]{gab} a localizing subcategory $\cL\subseteq \mathrm{Vect}^{\oplus I}$ is full on the objects
  \begin{equation}\label{eq:defq}
    \{x\in \mathrm{Vect}^{\oplus I}\ |\ \mathrm{hom}(x,Q)=0,\ \forall Q\}
  \end{equation}
  where $Q\in \mathrm{Vect}^{\oplus I}$ ranges over some family $\cF$ of injective objects. The defining condition for $x$ in \Cref{eq:defq} can be recast as requiring that the $i$-component of $x$ vanish for all $i\in I$ for which some $Q\in \cF$ has non-vanishing $i$-component. In short, $\cL$ must be of the form
  \begin{equation*}
    \mathrm{Vect}^{\oplus I'}\subseteq \mathrm{Vect}^{\oplus I}
  \end{equation*}
  for some subset $I'\subseteq I$ and hence we can identify
  \begin{equation*}
    \mathrm{Vect}^{\oplus I}/ker\pi \simeq \mathrm{Vect}^{\oplus (I\setminus I')}. 
  \end{equation*}
  To avoid overburdening the notation, we may as well assume (as we will for the duration of the proof) that
  \begin{equation*}
    \pi:\mathrm{Vect}^{\oplus I}\to \cC
  \end{equation*}
  is faithful. Together with $\pi\iota\cong \id_{\cC}$, that faithfulness then implies that
  \begin{equation}\label{eq:iota}
    \iota:\cC\to \mathrm{Vect}^{\oplus I}
  \end{equation}
  is {\it full} (in addition to being a faithful left adjoint). In other words, \Cref{eq:iota} identifies $\cC$ with a full subcategory of $\mathrm{Vect}^{\oplus I}$, closed under forming colimits.

  With all of this in hand, we define the {\it support} $\mathrm{supp}_{I}(\cC)$ of $\cC$ in $I$ to be the subset of those $i\in I$ for which there is some object $x\in \cC$ whose $i$-component (a vector space) is non-zero. It is easy to construct, for each $i\in \mathrm{supp}_{I}(\cC)$, an endomorphism $x\to x$ in $\cC$ whose cokernel is precisely the one-dimensional $i$-indexed vector space $k_i$ (regarded as an object of $\mathrm{Vect}^{\oplus I}$). But then, since \Cref{eq:iota} is
  \begin{itemize}
  \item full and
  \item colimit-closed,
  \end{itemize}
  $k_i\in \cC$. Taking colimits now shows that in fact $\cC$ (regarded as a full subcategory of $\mathrm{Vect}^{\oplus I}$ via \Cref{eq:iota}) is precisely
  \begin{equation*}
    \mathrm{Vect}^{\oplus\mathrm{supp}_I(\cC)}\subseteq \mathrm{Vect}^{\oplus I}. 
  \end{equation*}
  This finishes the proof, taking $S=\mathrm{supp}_I(\cC)$.
\end{proof}

In particular, taking $I$ to be a singleton in the proof of \Cref{th:naive}, we obtain the following characterization of simply-dualizable 2-$k$-modules.

\begin{corollary}\label{cor:simp}
  A non-trivial simply 1-dualizable 2-$k$-module $\cC$ is equivalent to $\mathrm{Vect}$ as a 2-$k$-module.
  \qedhere
\end{corollary}

\begin{example}\label{ex:notiso}
  In \Cref{eq:pipi} let
  \begin{itemize}
  \item $I=\{1,2\}$ be a two-elements set, so that $\mathrm{Vect}^{\oplus I}$ is equivalent to $\cM_{k\times k}$;
  \item $\cC\simeq \mathrm{vect}$;
  \item $\pi:\cM_{k\times k}\to \mathrm{Vect}$ be scalar restriction via the diagonal map $k\to k\times k$;
  \item $\iota:\mathrm{Vect}\to \cM_{k\times k}$ be scalar restriction via the projection $k\times k\to k$ on, say, the first component.
  \end{itemize}
  Then $\pi$ is clearly faithful (so that $T$ in \Cref{eq:pipi} is an equivalence), $\pi$ and $\iota$ are both left adjoints, we have $\pi\iota\cong \id_{\cC}$, but $\pi$ is not an equivalence.  
\end{example}

\section{Singleton categories}\label{se:sing}

In the context of \Cref{cor:simp} one considers 1-dualizable 2-$k$-modules where every object is a coproduct of copies of a single object $x$. The result concludes that $\cC$ is equivalent to $\mathrm{Vect}$ bu it seems pertinent, in view of that setup, to examine such ``singly-generated'' categories in more detail.

\begin{definition}\label{def:sing}
  A {\it singleton category} is a (usually linear) category $\cC$ with a distinguished object $x$ such that all objects of $\cC$ are (possibly empty) coproducts of copies of $x$. The object $x$ itself is then a {\it singleton} for $\cC$.

  The same term applies to narrower classes of categories: singleton 2-$k$-modules, singleton Grothendieck categories, etc.
\end{definition}

As seen above, being singleton is a weakening of the notion of simple 1-dualizability. First, recall the following notion (e.g. \cite[Chapter 4, Proposition 2.3]{pop}).

\begin{definition}\label{def:spec}
  A Grothendieck category is {\it spectral} if all of its objects are injective or, equivalently, projective. 
\end{definition}

Our first observation is

\begin{lemma}\label{le:spec}
A singleton Grothendieck category is spectral in the sense of \Cref{def:spec}. 
\end{lemma}
\begin{proof}
  Let $\cC$ be a singleton Grothendieck category and $x\in \cC$ an object with the property that every other object of $\cC$ is a copower $x^{\oplus S}$ (for some set $S$).
  
  Note that there are arbitrarily large injective copowers $x^{\oplus S}$ (i.e. for arbitrarily large sets $S$): simply take an injective envelope $E$ of some copower of $x$ whose set of subobjects has some arbitrarily large cardinality, and use the fact that by assumption, $E$ must itself be a copower of $x$.

  Finally, embed an arbitrary copower $x^{\oplus S}$ as a summand into an injective $x^{\oplus T}$, $|S|\le |T|$, hence concluding that the arbitrary copower $x^{\oplus S}$ must again be injective. 
\end{proof}

\Cref{le:spec} is what first goes wrong as soon as one tries to generalize these considerations to more than a single object $x$:

\begin{example}
  Let $k$ be a field. The category $\cM_A$ of modules over the ring $A=k[\varepsilon]/(\varepsilon^2)$ of dual numbers has the property that every object is a direct sum of copies of $k$ and $A$ (the indecomposable objects of $\cM_A$). $\cM_A$ is not spectral, showing that even two ``base'' objects will undo \Cref{le:spec}.
\end{example}

The following remark justifies the phrase ``{\it the} singleton''.

\begin{lemma}\label{le:singuniq}
  In a singleton Grothendieck category all singletons are isomorphic. 
\end{lemma}
\begin{proof}
  If $x$ and $y$ are singletons, each is a summand of the other, and they are moreover injective by \Cref{le:spec}. It follows that they are isomorphic, e.g. by \cite[Theorem 10.14]{good} (which is phrased for modules, but goes through in arbitrary Grothendieck categories).
\end{proof}

\begin{lemma}\label{le:xsing}
  Let $\cC$ be a singleton Grothendieck category and $0\ne x\in \cC$ a singleton. Then, all non-zero subobjects of $x$ are isomorphic. 
\end{lemma}
\begin{proof}
  The argument is similar to that used in the proof of \Cref{le:singuniq}: if $y\le x$ is a subobject then $x$ and $y$ are injectives realizable as summands in each other, and hence are isomorphic by \cite[Theorem 10.14]{good}.
\end{proof}

Now let $\cC$ be a singleton Grothendieck category. \Cref{le:spec} imposes strong restrictions on the structure of $\cC$. Specifically, \cite[Chapter 4, Theorem 2.8]{pop} says that there is a regular self-injective ring $A$ such that
\begin{equation}\label{eq:cmar}
  \cC\simeq \cM_A^r,
\end{equation}
where the latter (with `$r$' standing for `reduced') denoting the full subcategory
\begin{equation*}
  \cM_A^r\subset \cM_A
\end{equation*}
consisting of right modules appearing as direct summands of powers $A^I$ for sets $I$ (see also \cite[Chapter XII, Theorem 1.3]{stns}). Furthermore, the equivalence \Cref{eq:cmar} can be described concretely as
\begin{equation}\label{eq:xcmar}
  \mathrm{Hom}_{\cC}(x,-):\cC\to \cM_A^r
\end{equation}
for some generator $x\in \cC$.

Regarded as a functor to the category $\cM_A$ of {\it all} modules $\mathrm{Hom}_{\cC}(x,-)$ is a right adjoint and hence preserves products. It does {\it not}, in general, preserve coproducts: although the category $\cM_A^r$ is cocomplete (indeed, it is Grothendieck), direct sums are somewhat more involved that direct products, which are simply the usual Cartesian products of modules. The following result follows, for instance, from \cite[Proposition 9.1]{good}, and we omit the proof.

\begin{lemma}\label{le:sums}
  For a regular self-injective ring $A$, the coproduct in $\cM_A^r$ of a family of objects $X_j\in \cM_A^r$, $j\in J$ is the unique summand $E$ of the product $\prod_j X_j$ fitting into the chain
  \begin{equation*}
    \bigoplus_j X_j \le E \le_{\oplus} \prod_j X_j
  \end{equation*}
  so that the left hand inclusion is essential.   
  \qedhere
\end{lemma}

\begin{remark}
  In particular, the uniqueness of such a summand is part of the statement; it is that uniqueness that follows from the cited result, specifically \cite[Proposition 9.1 (e)]{good}.
\end{remark}

Throughout the discussion we fix a singleton Grothendieck category $\cC$ with a singleton $0\ne x\in \cC$. $A$ typically denotes a regular right self-injective ring with $\cC\simeq \cM_A^r$ but we do {\it not} assume, in general, that the equivalence identifies $A$ to $x$ (or to a singleton, in general).

The goal is to understand the structural restrictions imposed on $A$ by the demand that $\cC$ be singleton. A first observation, \Cref{le:fact} below, borrows the following operator-algebraic term (e.g. \cite[Chapter 1, Exercise 7]{kapl-op}, \cite[Definition 3.1.2]{mv1}, \cite[Definition II.3.2]{tak1}, etc.).

\begin{definition}\label{def:fact}
  A regular ring is a {\it factor} if its only central idempotents are $0$ and $1$. 
\end{definition}

\begin{lemma}\label{le:fact}
The regular self-injective ring $A$ in \Cref{eq:cmar} is a factor. 
\end{lemma}
\begin{proof}
  Let $e\in A$ be a central idempotent. If $e$ were non-trivial (i.e. neither $0$ nor $1$) then $eA$ and $(1-e)A$ would both be objects of $\cM_A^r\simeq \cC$ whose annihilators ($(1-e)A$ and $eA$ respectively) are not both contained in the annihilator of a non-trivial object in $\cM_A^r$.

  On the other hand, $x\in \cC$ is a summand of every non-zero object in $\cC$ and hence, as an $A$-module, its annihilator contains those of every other $\cC$-object. This gives the desired contradiction, proving that $e\in \{0,1\}$.
\end{proof}

In particular,

\begin{corollary}\label{cor:centfield}
In the context of \Cref{le:fact} the center of $A$ is a field.   
\end{corollary}
\begin{proof}
  This follows from \Cref{le:fact} and \cite[Corollary 1.15]{good}.
\end{proof}

In a sense, \Cref{le:xsing,le:fact} completely classify singleton Grothendieck categories. First, we need

\begin{definition}\label{def:rngsing}
  A right $A$-module is {\it isominimal} if it is isomorphic to all of its non-zero summands. The ring $A$ itself is isominimal if it is isominimal as a right $A$-module.

Finally, $x\in A$ is isominimal if $xA$ is.
\end{definition}

\begin{remark}
  Note that isominimal regular rings are automatically factors. 
\end{remark}

\begin{theorem}\label{th:clssing}
  For a non-trivial Grothendieck category $\cC$, the following conditions are equivalent:
  \begin{enumerate}[(1)]
  \item\label{item:6} $\cC$ is singleton.
  \item\label{item:7} $\cC$ is equivalent to $\cM_A^r$ for a (non-zero) regular isominimal self-injective ring.  
  \item\label{item:8} $\cC$ is equivalent to $\cM_A^r$ for a regular self-injective factor containing some non-zero isominimal idempotent $e\in A$.
  \end{enumerate}
\end{theorem}
\begin{proof}

  {\bf \Cref{item:6} $\Rightarrow$ \Cref{item:7}.} This is a consequence of \Cref{le:xsing} and the surrounding discussion: we have an equivalence \Cref{eq:xcmar} for a regular self-injective $A$, identifying a singleton $x\in \cC$ to $A\in \cM_A^r$. It follows from \Cref{le:xsing} that all non-zero principal right ideals of $A$ are isomorphic to $A$, and thus the latter must be isominimal.

{\bf \Cref{item:7} $\Rightarrow$ \Cref{item:8}.} Trivial.

{\bf \Cref{item:8} $\Rightarrow$ \Cref{item:6}.} We will argue that every summand
\begin{equation*}
  X\le_{\oplus} A^I
\end{equation*}
of a power of $A$ is a direct sum (in $\cM_A^r$) of copies of $eA$. We first do this for the one-element set $I$:

{\bf Claim: Non-zero summands of $A$ are direct sums in $\cM_A^r$ of copies of $eA$.} Let $xA\le_{\oplus A}$ be a non-zero summand 
\begin{equation}\label{eq:xjmax}
  X_j\le xA,\ j\in J
\end{equation}
be a maximal independent family of objects $X_j\cong eA$ (`independent' in the sense that, as before, their sum in $A$ is direct). Given the description of coproducts in $\cM_A^r$ provided by \Cref{le:sums}, this means proving that the actual, module-theoretic direct sum $\bigoplus X_j$ is essential in $xA$.

If not, the unique summand $E\le_{\oplus}xA$ that contains $\bigoplus_j X_j$ essentially has a non-zero complement $F$, i.e. $xA=E\oplus F$. Since $A$ is a regular self-injective factor, it satisfies the comparability axiom (\cite[p.80, Definition]{good}) by \cite[Corollary 9.16]{good} and hence one of $eA$ and $F$ must be a summand of the other. Either way, the isominimality of $e$ implies that $eA$ is realizable as a summand of $F$, in which case the independence of the family
\begin{equation*}
  \{X_j,\ j\in J\}\sqcup \{F\}
\end{equation*}
contradicts the maximality of \Cref{eq:xjmax}. This concludes the proof of the claim.

We now resume the main line of reasoning. Once more, consider a maximal independent family $X_j\cong eA$ of copies of $eA$. We claim that
\begin{equation*}
  X=\bigoplus X_j\text{ in }\cM_A^r. 
\end{equation*}
As above, we have to show that the module-theoretic direct sum $\bigoplus X_j$ is essential in $X$. The argument is entirely parallel: assuming the opposite, the unique summand $E\le_{\oplus}X$ that contains $\bigoplus_j X_j$ essentially will be proper in $X$ and hence have a complement:
\begin{equation*}
  X=E\oplus F,\quad F\ne 0.
\end{equation*}
$F$ embeds as a summand in the original power $A^I$, and hence maps non-trivially to some factor $A^I\to A$ via a morphism $\varphi:F\to A$. We know from \cite[Proposition 9.1, parts (a) and (b)]{good} that
\begin{equation*}
  0\to \ker~\varphi\to F\to \mathrm{im}~\varphi\to 0
\end{equation*}
splits and that the image is a summand of the codomain $A$ of $\varphi$. It thus follows that a non-zero summand $Z$ of $A$ can be realized as a summand of the complement $F=X\ominus E$. But this means that the maximal family $\{X_j\}$ can be extended by some summand $eA\le_{\oplus} Z$ (at least one exists by the claim proven above). This contradicts the maximality of the family $\{X_j\}$ and thus finishes the proof.
\end{proof}

As for a more concrete description of singleton categories, we can resort (via \Cref{eq:cmar}) to the classification theory of regular rings as covered in \cite[Chapter 10]{good}. 

\begin{proposition}\label{pr:13}
Let $\cC\cong \cM_A^r$ be a singleton category for a regular self-injective ring $A$ as the singleton. Then, $A$ is of type I or III. 
\end{proposition}
\begin{proof}
  Given that
  \begin{itemize}
  \item $A$ is a factor by \Cref{le:fact} on the one hand and 
  \item being regular self-injective, $A$ decomposes uniquely as a product of type-I, II and III rings on the other (\cite[Theorem 10.13]{good}), 
  \end{itemize}
  all we have to do is rule out type II.

  If $A$ were type-II it would have a non-zero directly finite idempotent $e\in A$ \cite[Proposition 10.8]{good}, in which case $eA$ could not be isomorphic to any of its proper direct summands by the very definition of direct finiteness for modules (\cite[p.49, Definition]{good}). This means that $eA$ must be a {\it finite} direct sum of copies of the singleton $x$, so we may as well assume $eA\cong x$.

  By direct finiteness again $x\cong eA$ cannot be isomorphic to a proper direct summand, so in fact $eA$ must be indecomposable. But then $e\ne 0$ is an abelian idempotent in the sense of \cite[p.110, Definition]{good}: the regular ring $eAe$ has no non-trivial idempotents at all, and hence is a division ring (and in particular an abelian regular ring). This contradicts the assumption that $A$ is type-II (\cite[p.113, Definition]{good} requires that no non-zero abelian idempotents exist), finishing the proof.
\end{proof}

As expected, the ``easy'' case is that of type I. 

\begin{proposition}\label{pr:divrng}
  For a Grothendieck category $\cC$ the following conditions are equivalent:
  \begin{enumerate}[(1)]
  \item\label{item:1} $\cC$ is equivalent to the category of right vector spaces over a division ring $D$.
  \item\label{item:2} $\cC$ is a singleton Grothendieck category with an indecomposable singleton object.
  \item\label{item:3} $\cC$ is singleton and any regular self-injective ring $A$ acting as a singleton in $\cM_A^r\simeq \cC$ is a division ring.
  \item\label{item:4} $\cC$ is singleton and any regular self-injective ring $A$ with $\cM_A^r\simeq \cC$ is a factor of type I.
  \item\label{item:5} $\cC$ is equivalent to $\cM_B^r$ for {\it some} regular self-injective type-I factor $B$. 
  \end{enumerate}
\end{proposition}
\begin{proof}
  {\bf \Cref{item:1} $\Rightarrow$ \Cref{item:2}.} Indeed, in that case the (unique, up to isomorphism) singleton is $D$ itself, which as indecomposable as a right $D$-module.

  {\bf \Cref{item:2} $\Rightarrow$ \Cref{item:3}.} If there is an indecomposable singleton (which will automatically be unique up to isomorphism) then $A\in \cM_A$ must be indecomposable and hence simple as a right $A$-module (because principal right ideals are summands). It follows that $A$ is a division ring.

  {\bf \Cref{item:3} $\Rightarrow$ \Cref{item:4}.} Division rings are type-I.

  {\bf \Cref{item:4} $\Rightarrow$ \Cref{item:5}.} Trivial. 

  {\bf \Cref{item:5} $\Rightarrow$ \Cref{item:1}.} Being of type I, \cite[Theorem 10.6]{good} implies that $B$ is of the form $\mathrm{End}_D(M)$ where
  \begin{itemize}
  \item $D$ is an abelian regular self-injective ring (`abelian' in the sense of \cite[Definitions on p.25 and p.110]{good}: idempotents are central); 
  \item $M$ is a non-singular injective right $D$-module. 
  \end{itemize}
  Since $B$ is a factor so is $D$. Abelian factors are division rings, so $B$ is a full endomorphism rings of a vector space $M$ over a division ring $D$. But in that case we have
  \begin{equation*}
    \cC\simeq \cM_{\mathrm{End}_D(M)}^r\simeq \mathrm{Vect}_D,
  \end{equation*}
  finishing the proof.
\end{proof}

As for the type-III arm of \Cref{pr:13}, not only does it occur, but in fact accounts for everything \Cref{pr:divrng} doesn't.

\begin{proposition}\label{pr:iii}
  For any type-III regular self-injective factor $A$, the category $\cM_A^r$ is singleton. 
\end{proposition}
\begin{proof}
  \cite[Theorem 5-3.6]{gw} implies that for a type-III regular self-injective factor $A$ the poset formed by the isomorphism classes of finitely-generated $A$-projectives (with order given by being a direct summand) is isomorphic to an initial segment of ordinal numbers.

  In particular, such a factor has a non-zero isominimal idempotent in the sense of \Cref{def:rngsing} and the conclusion follows from \Cref{th:clssing}.
\end{proof}

Collecting together \Cref{pr:13,pr:divrng,pr:iii}, we obtain the following complete classification of singleton Grothendieck categories.

\begin{theorem}\label{th:complete}
  The singleton Grothendieck categories $\cC$ are precisely those of the form
  \begin{enumerate}[(a)]
  \item $\cM_A^r\simeq \mathrm{Vect}_D$ for a division ring $A=D$;
  \item $\cM_A^r$ for a type-III regular self-injective factor or, equivalently, simple ring $A$.
  \end{enumerate}
\end{theorem}
\begin{proof}
  Given the preceding discussion, the only thing left to prove is that in the type-III case the ring $A$ can always be chosen to be simple. To see this, recall first that we can identify $\cC$ with $\cM_A^r$ via \Cref{eq:xcmar}, thus taking for our $A$ the endomorphism ring $\mathrm{End}_{\cC}(x)$ of a singleton. It now remains to observe that endomorphism rings of singletons are simple:

  According to \Cref{le:xsing} all non-zero subobjects of $x$ are mutually isomorphic (and hence isomorphic to $x$). This means that a non-zero endomorphism $f:x\to x$ implements an isomorphism of a summand $y$ of $x$ onto another summand $z$, with
  \begin{equation*}
    x\cong y\cong z. 
  \end{equation*}
  But then, composing $f$ on either side with
  \begin{itemize}
  \item a morphism $\psi:x\cong y\subseteq x$ and similarly,
  \item a morphism $\theta:x\to x$ that maps $z\le x$ isomorphically onto $x$
  \end{itemize}
  we obtain $\id_x = \theta\circ f\circ \psi$, hence the conclusion that $\mathrm{End}_{\cC}(x)$ is simple.
\end{proof}



\begin{thebibliography}{10}

\bibitem{ar}
Ji\v{r}\'{\i} Ad\'{a}mek and Ji\v{r}\'{\i} Rosick\'{y}.
\newblock {\em Locally presentable and accessible categories}, volume 189 of
  {\em London Mathematical Society Lecture Note Series}.
\newblock Cambridge University Press, Cambridge, 1994.

\bibitem{ah}
Ulrich Albrecht and Jutta Hausen.
\newblock Nonsingular modules and {$R$}-homogeneous maps.
\newblock {\em Proc. Amer. Math. Soc.}, 123(8):2381--2389, 1995.

\bibitem{bd}
John~C. Baez and James Dolan.
\newblock Higher-dimensional algebra and topological quantum field theory.
\newblock {\em J. Math. Phys.}, 36(11):6073--6105, 1995.

\bibitem{bcjf}
Martin Brandenburg, Alexandru Chirvasitu, and Theo Johnson-Freyd.
\newblock Reflexivity and dualizability in categorified linear algebra.
\newblock {\em Theory Appl. Categ.}, 30:Paper No. 23, 808--835, 2015.

\bibitem{cjf}
Alex Chirvasitu and Theo Johnson-Freyd.
\newblock The fundamental pro-groupoid of an affine 2-scheme.
\newblock {\em Appl. Categ. Structures}, 21(5):469--522, 2013.

\bibitem{gab}
Pierre Gabriel.
\newblock Des cat\'{e}gories ab\'{e}liennes.
\newblock {\em Bull. Soc. Math. France}, 90:323--448, 1962.

\bibitem{good}
K.~R. Goodearl.
\newblock {\em von {N}eumann regular rings}.
\newblock Robert E. Krieger Publishing Co., Inc., Malabar, FL, second edition,
  1991.

\bibitem{gw}
K.~R. Goodearl and F.~Wehrung.
\newblock The complete dimension theory of partially ordered systems with
  equivalence and orthogonality.
\newblock {\em Mem. Amer. Math. Soc.}, 176(831):vii+117, 2005.

\bibitem{kapl-op}
Irving Kaplansky.
\newblock {\em Rings of operators}.
\newblock W. A. Benjamin, Inc., New York-Amsterdam, 1968.

\bibitem{kel}
G.~M. Kelly.
\newblock Basic concepts of enriched category theory.
\newblock {\em Repr. Theory Appl. Categ.}, (10):vi+137, 2005.
\newblock Reprint of the 1982 original [Cambridge Univ. Press, Cambridge;
  MR0651714].

\bibitem{lur}
Jacob Lurie.
\newblock On the classification of topological field theories.
\newblock In {\em Current developments in mathematics, 2008}, pages 129--280.
  Int. Press, Somerville, MA, 2009.

\bibitem{macl}
Saunders Mac~Lane.
\newblock {\em Categories for the working mathematician}, volume~5 of {\em
  Graduate Texts in Mathematics}.
\newblock Springer-Verlag, New York, second edition, 1998.

\bibitem{mv1}
F.~J. Murray and J.~Von~Neumann.
\newblock On rings of operators.
\newblock {\em Ann. of Math. (2)}, 37(1):116--229, 1936.

\bibitem{pop}
N.~Popescu.
\newblock {\em Abelian categories with applications to rings and modules}.
\newblock Academic Press, London-New York, 1973.
\newblock London Mathematical Society Monographs, No. 3.

\bibitem{stns}
Bo~Stenstr\"{o}m.
\newblock {\em Rings of quotients}.
\newblock Springer-Verlag, New York-Heidelberg, 1975.
\newblock Die Grundlehren der Mathematischen Wissenschaften, Band 217, An
  introduction to methods of ring theory.

\bibitem{tak1}
M.~Takesaki.
\newblock {\em Theory of operator algebras. {I}}, volume 124 of {\em
  Encyclopaedia of Mathematical Sciences}.
\newblock Springer-Verlag, Berlin, 2002.
\newblock Reprint of the first (1979) edition, Operator Algebras and
  Non-commutative Geometry, 5.

\end{thebibliography}

\addcontentsline{toc}{section}{References}

\Addresses

\end{document}